\documentclass[11pt]{article}%
\usepackage{enumerate}
\usepackage{enumitem}
\usepackage{amsmath,amssymb,amsfonts}
\usepackage{amsfonts}%
\usepackage{graphicx}
\usepackage{color}
\usepackage{pdfsync}

\usepackage{amssymb}
\usepackage{amsmath,amsfonts,amsthm}

\usepackage[latin1]{inputenc}

\newtheorem{thm}{Theorem}[section]

\newtheorem{lem}{Lemma}[section]

\newtheorem{rmk}{Remark}[section]

\let\Section=\section

\def\section{\setcounter{equation}{0}\Section}

\def\nd{\noindent}

\def\hsf{\hspace*{\fill}}

\def\proof{{\rm \bf Proof}}

\def\r{\bf{R}}

\begin{document}
\title{Multiple Sign Changing Radially Symmetric Solutions\\
 in  a General Class  of Quasilinear Elliptic Equations \thanks{Partially supported by PROCAD/CAPES/UFG/UnB-Brazil}}

\author{Claudianor O. Alves\thanks{ Claudianor O. Alves was  supported  by CNPQ/Brasil. }~~ Jose V. A. Goncalves\thanks{Jose V. A.  Goncalves was  supported by CNPq/Brazil}~~  Kaye O. Silva\thanks{Kaye O. Silva was  supported by CAPES/Brazil}}

\date{}

\pretolerance10000

\maketitle

\begin{abstract}
\noindent {\small  
In this paper we prove that the equation $ -( r^\alpha\phi(|u'(r)|)u'(r))' = \lambda r^\gamma f(u(r)), ~0<r<R$, where $\alpha, \gamma, {\bf{R}}$ are given real numbers,  $\phi : (0, \infty) \to (0, \infty)$ is a  suitable twice  differentiable function, $\lambda > 0$ is a real parameter and  $f:{\bf{R}}\to{\bf{R}}$ is continuous,  admits an infinite sequence of sign-changing solutions satisfying  $u'(0) =u(R) =0$. The function $f$ is  required to satisfy $tf(t)>0$ for  $ t\neq 0$. Our technique  explores fixed point arguments applied to suitable integral equations and shooting arguments. Our main result extends earlier ones in the case $\phi$ is in the form $\phi(t) = |t|^{\beta}$ for an appropriate constant $\gamma$. }
\end{abstract}

\section{Introduction}

\nd  We study the nonlinear eigenvalue problem

~\begin{equation}\label{P}\tag{$P_\lambda$}
\left\{\begin{array}{rllr}
-( r^\alpha\phi(|u'(r)|)u'(r))' &=& \lambda r^\gamma f(u(r)), ~~0<r<R,\\
\\
u'(0) &=&u(R) =0,
 \end{array}
\right.
\end{equation}

\nd where  $\lambda > 0$ is a parameter, $f:\mathbb{R}\to\mathbb{R}$ is continuous and $\alpha,\gamma\in\mathbb{R}$  are suitable constants.
\vskip.2cm

\nd We shall assume that  $\phi : (0, \infty) \to (0, \infty)$ is a  twice  differentiable, $C^1$-function,  satisfying

\begin{itemize}
  \item[($\phi_1$)]  \     \ $\mbox{(i)} \ \  t \phi(t) \to 0 \ \mbox{as} \  t\to 0, \\
\\
 ~~~ \mbox{(ii)} \  t \phi(t) \to\infty \ \mbox{as} \  t\to\infty$,
  \item[($\phi_2$)]  \ $t \phi(t) \ \mbox{ is strictly increasing in}~ (0, \infty)$,
  \item[($\phi_3$)] there  are constants  $\gamma_1, \gamma_2>1$ such that $$\gamma_1-1\le\frac{(t \phi(t))^{\prime}}{\phi(t)}\le\gamma_2-1,\ t > 0.$$
 \end{itemize}

\nd Concerning $f$, the following conditions will be imposed:

 \begin{itemize}
  \item[($f_1$)] $tf(t)>0,\ t\neq 0$,
  \item[($f_2$)]  there exists $d_\infty>0$ such that $f$ is nondecreasing in $(-\infty,d_\infty]$,
  \item[($f_3$)]~~~~~  $\displaystyle \liminf_{\nu\to 0^{\pm}}\int_0^\nu |f(t)|^{\frac{-1}{\gamma_1-1}}\operatorname{sgn} (f(t))dt< \infty$.
 \end{itemize}

\begin{rmk} We observe that condition ($f_3$) is equivalent to the following:
\begin{itemize}
  \item[($f_3'$)]~~~~~~~~~~~ $\max \Big \{ \int_{-x}^0 [-f(t)]^{\frac{-1}{\gamma_1-1}}dt,~ \int_{0}^y [f(t)]^{\frac{-1}{\gamma_1-1}}dt \Big \}<\infty $,
\end{itemize}
 \nd  $\mbox{for each} \ x,y>0$, where  $\gamma_1'=\gamma_1/(\gamma_1-1)$.
\end{rmk}

\vskip.3cm
 \nd Our main objective in this work is to prove the following  result:

 \begin{thm}\label{T1} Let $f\in C(\mathbb{R})\cap C^1(\mathbb{R}\setminus\{0\})$. Assume ($\phi_1$)-($\phi_3$),  ($f_1$)-($f_2$) and
\begin{equation}\label{alphabeta}
\tag{$\gamma,\alpha$} \gamma\ge \max\left\{\alpha,\frac{-\alpha}{\gamma_1-1}\right\}.
\end{equation}

\nd Then there is a positive number $\Lambda$ such that  for each $\lambda\in (0,\Lambda]$,
problem (\ref{P}) admits a positive solution $u_0$ and an infinite  sequence $\{u_\ell\}_{\ell=1}^\infty$ of solutions satisfying:
\begin{equation}\label{T12}
u_l(0) =d_\ell,
\end{equation}
\begin{equation}\label{T13}
u_\ell\ \mbox{has precisely $\ell$ zeroes in } (0,R),
\end{equation}

\nd where $\{d_\ell\}_{\ell=1}^\infty$ is an infinite  sequence of real numbers  such that

 \begin{equation}\label{T11}   d_\infty > d_1>\cdots>d_\ell>\cdots>0.\end{equation}
\end{thm}
\vskip.2cm

 \nd  The proof of  Theorem \ref{T1}  is strongly based on  the shooting method. In this regard,  consider the initial value problem
\begin{equation}\label{PP}\tag{$P_{\lambda,d}$}
\left\{\begin{array}{rllr}
-( r^\alpha\phi(|u'(r)|)u'(r))' &=& \lambda r^\gamma f(u(r)),~~ ~~r>0,\\
\\
u(0)=d,~~ u'(0)=0,
 \end{array}
\right.
\end{equation}

\nd where  $d\in (0,d_\infty]$.
\vskip.3cm

\nd   The auxiliary result below will play a crucial role in this work. 

\begin{thm}\label{T2} Let $f\in C(\mathbb{R})\cap C^1(\mathbb{R}\setminus\{0\})$. Assume ($\phi_1$)-($\phi_3$), (\ref{alphabeta}) and ($f_1$)-($f_2$).  Then there exists a positive number $\Lambda= \Lambda(d_\infty) $ such that for each $\lambda\in (0,\Lambda]$, problem (\ref{PP}) has a unique solituon $u(\cdot,d,\lambda)=u(\cdot,d) \in C^1([0,\infty))$. In addition, for each $d \in (0,d_\infty]$,  there is a sequence  $\{ z_{\ell} \}_{{\ell=1}}^{{\infty}}$ of zeroes of  $u(\cdot,d)$,   $z_{\ell}=z_{\ell}$(d), such that

 \begin{equation}\label{T23}
\begin{array}{l}
z_1(d_\infty) \geq R,~~u(r,d)>0\ \mbox{if}\ 0<r<z_1(d),\\
\\
 z_1(d) < z_2(d) < \cdots < z_\ell(d) < \cdots,\\
\\
 u'(r,d)<0\ \mbox{if}\ 0<r\le z_1(d), u(r,d)\neq 0\ \mbox{if}\ z_{\ell}<r<z_{\ell+1}~ \mbox{and}~ u'(z_{\ell},d) \neq 0, \\
 \end{array}
\end{equation}
\vskip.2cm

\begin{equation}\label{T24}
\begin{array}{l}
z_{\ell}(d)\to 0\ \mbox{as}\ d\to 0\ \mbox{and}\ z_{\ell}(d)\to z_{\ell}( {\underline{d}})\ \mbox{as}\ d\to {\underline{d}},~~ {\underline{d}} \in (0,d_\infty],
 \end{array}
\end{equation}

\begin{equation}\label{T25}\begin{split}&\mbox{if}\ {\underline{d}} \in (0,d_\infty] \ \mbox{and}\ u(\cdot, {\underline{d}})\ \mbox{has}\ k\ \mbox{zeroes in}\ (0,R)\ \mbox{then}\ u(\cdot,d)\ \mbox{has at most}\ k+1 \\ & \mbox{zeroes in}\ (0,R)\ \mbox{whenever}\ d< {\underline{d}}, d\ \mbox{close enough to}\  {\underline{d}}.\end{split}\end{equation}
\end{thm}

\section{Background}

\nd Consider the problem
\begin{equation}\label{P1}\tag{$P1_\lambda$}
\left\{\begin{array}{rllr}
-{\rm div}(a(x) \vert \nabla{u(x)} \vert^{\beta} \nabla{u(x)} ) & =&  \lambda~ b(x) f(u),~~ x \in B_R,\\
\\
u(x) &=& 0,~~ x \in \partial{B_R},
\end{array}
\right.
\end{equation}
\nd where $B_{R} \subset {\r^N}$ is the ball of radius $R$ centered at the origin, the functions  $a,b$ are radially symmetric and  $\beta > - 1$. Making $a = b \equiv 1$, $\beta = p - 2$ with $1 < p < \infty$ and $\lambda = 1$, (\ref{P1})  becomes
\begin{equation}\label{P2}\tag{$P2$}
\left\{\begin{array}{rllr}
\displaystyle - (r^{N - 1} |u'(r)|^{p - 2}u'(r))^{'}  &=&  r^{N - 1} f(u(r)),~~  0 <  r < R,\\
\\
 u'(0)  =  u(R) = 0.
\end{array}
\right.
\end{equation}

\nd It was shown in   \cite{JIa}    that  if $f(t) = |t|^{\delta-1}t$ with $1 < \delta + 1 < p < N$ then (\ref{P2}) has infinitely many nodal solutions.
\vskip.2cm

\nd In   \cite{mg}, it was shown that  the more general problem
\begin{equation}\label{P3}\tag{$P3$}
\left\{\begin{array}{rllr}
\displaystyle - (r^{\alpha} |u'(r)|^{\beta}u'(r))^{'}  &=&  \lambda  r^{\gamma} f(u(r)),~~  0 <  r < R,\\
\\
 u'(0)  =  u(R) = 0
\end{array}
\right.
\end{equation}

\nd admits infinitely many solutions if $\lambda$ is positive and small enough,

\begin{equation}
\beta > -1,~~ \gamma\ge \max\left\{\alpha,\frac{-\alpha}{\beta + 1}\right\},
\end{equation}

\nd and conditions $(f_1),(f_2)$ and a stronger form of $(f_3)$ hold.
\vskip.3cm

\nd Regarding (\ref{P3}), an example of a function safistying $(f_1),(f_2),(f_3)$  with  $\beta >0$ is $f(t) = arctg(t)$.
\vskip.3cm

\nd As was pointed out by Clement, Figueiredo \& Mitidieri \cite{CFM} the  operator
$$
(r^{\alpha} |u'(r)|^{\beta}u'(r))^{'}
$$
\nd represents the radial form of  the well known operators:
\vskip.3cm
\begin{description}
\item{} {\it p-Laplacian  with $1<p<N$~when~  $\alpha=N-1,~  \beta= p-2$,}
\item{} {\it  k-Hessian  with $1\leq k \leq N$~ when~ $ \alpha= N-k,~ \beta= k-1.$}
\end{description}

\nd The problem   
\begin{equation}\label{Phi}\tag{$\Phi$}
\left\{\begin{array}{rllr}
-\Delta_{\Phi} u &=& \lambda f(u)~~ \mbox{in}~ B\\
\\
u &=&  0~~ \mbox{on}~   \partial{B},
\end{array}
\right.
\end{equation}

\nd where
$$
\Phi(t)=\int_0^t s\phi(s)ds,
$$
\nd  $\Delta_{\Phi}$ is the  $\Phi$-Laplacian operator namely
$$
\Delta_{\Phi} u = {\rm div} \big(\phi(\mid\nabla u\mid)\nabla u\big),
$$
\nd and $B \subset {R}^N$ is the ball of radius $R$  centered at the origin, was addressed by many authors (see e.g. Fukagai   \&  Narukawa\cite{fn-1} and its references).
A weak  solution of  ($\Phi$) is by definition an element $u \in W_0^{1, \Phi}({B})$  (the usual Orlicz-Sobolev space) such that
\begin{equation}\label{radialphi}
\int_B \phi(|\nabla u|) \nabla u \nabla v  dx =\lambda\int_B f(u) v dx,~  v\in W_0^{1, \Phi}({B}).
\end{equation}

\nd The radially symmetric form of $(\Phi)$ is
$$
\left\{\begin{array}{rllr}
\displaystyle - (r^{N-1}\phi(|u'(r)|)u'(r))'  &=&  \lambda r^{N-1}f(u(r),~ 0 <  r < R\\
\\
 u'(0)  =  u(R) = 0
\end{array}
\right.
$$
\nd which is a special case of (\ref{P}), (see further  remarks  in the Appendix).
\vskip.2cm

 \nd  Theorem 1.1 extends the main results of  {\rm \cite{mg}, \cite{JIa}},   in the sense that we were able to treat both with a more general class of quasilinear operators  and a broader class of terms $f$.
\vskip.2cm

\nd Problems like (\ref{P}) have been investigated by many authors and we would like to refer the reader to
Saxton  $\&$  Wei {\rm \cite{SWei}}, Castro $\&$  Kurepa  {\rm \cite{CKr}}, Cheng {\cite{Cheng}}, Strauss {\rm \cite{St}}, Ni $\&$ Serrin  {\rm \cite{NiSer}},  Castro, C\'ossio $\&$  Neuberger {\rm \cite{CCN}},  Fukagai  \& Narukawa \cite{fn-1}, Mihailescu \& Radulescu \cite{MR1,MR2} and their references. Here, we would point out that in \cite{fn-1}, Fukagai  \& Narukawa have mentioned that this type of problem appears in some fields of physics, such as, nonlinear elasticity, plasticity and generalized Newtonian fluids.

\section{Proof of Theorem \ref{T1}}

\nd Take $\lambda\in (0,\Lambda]$ where $\Lambda>0$ is given in Theorem \ref{T2}. We proceed in two steps:
\vskip.3cm

\nd {\bf Step 1}. {\bf {(Existence of a positive solution of (\ref{P}).)}} Let $d \in (0, d_\infty]$. We shall use the notations in Theorem  \ref{T2}. So $z_1 = z_1(d)$ denotes the first zero of $u(\cdot) = u(\cdot,d)$.
 Set 
$$
A_0=\Big \{d\in (0,d_\infty]~|~  z_1(d)\ge R \Big\}~\mbox{and}~ d_0 :=\inf A_0.
$$ 

\nd By (\ref{T23}) in  theorem \ref{T2},    $z_1(d_\infty)\ge R$.  So $A_0\neq \emptyset$. We will show that

\begin{equation}\label{th11}
d_0 > 0~\mbox{and}~ z_1(d_0)=R.
\end{equation}

\nd Indeed,  assume  on the contrary that  $d_0=0$. Take a sequence  $(d_j) \subseteq  A_0$ such that $d_j\to 0$.  By  (\ref{T24}), $z_1(d_j)\to 0$, which is a contradiction. 
\vskip.2cm

\nd Now,  assume $z_1(d_0)>R$.  Pick a sequence $(d_j) \subseteq  (0, d_{\infty}]$ such that  $d_j < d_0$ and $d_j\to d_0$. Applying  (\ref{T24}) we infer that   $z_1(d_j)\to z_1(d_0)$. Once $z_1(d_0)>R$, it follows that $z_1(d_j)>R$, which shows that $d_j \in A_{0}$. But this contradicts the definition of $d_0$.   Therefore $z_1(d_0)=R$ and this completes the proof of  (\ref{th11}). As a byproduct there is a positive solution of  (\ref{P}).
\vskip.3cm

\nd {\bf Step 2}.   {\bf {(Existence of an infinite sequence of sign-changing  solutions of (\ref{P}).)}} At first consider
$$
A_1 := \Big \{d\in (0,d_0]~|~ \ z_1(d)<R,\ z_{2}(d)\ge R \Big \}~~\mbox{and}~~ d_1 := \inf A_1.
$$ 
\nd We claim  that

\begin{equation}\label{th12}\begin{split}
A_1 \neq \emptyset,~~0 < d_1 <  d_0,\\
\\ 
\ z_1(d_1)<R,~~ z_{2}(d_1)=R.
\end{split}\end{equation}

\nd Let us show at first that $A_1\neq\emptyset$. Indeed, by {\bf Step 1} $z_1(d_0) = R$. By (\ref{T25}), if $d \in(0,d_0)$ with $d$ close to $d_0$ then $u(\cdot,d)$ has at most one zero in $(0,R)$.  Assume by contradiction that $u(\cdot,d)$ has no zero in $(0,R)$. Then $z_1(d) \geq R$ with $d < d_0$, impossible. It follows that $u(\cdot,d)$ has precisely one zero in $(0,R)$ and so $d \in A_1$, showing that $A_1\ne\emptyset$. 
\vskip.1cm

\nd To show that  $d_1>0$, assume by contradiction that there is a  sequence $\{d_j \} \subset  A_1$ such that  $d_j\to 0$. By (\ref{T24}), $z_2(d_j)\to 0$ contradicting $z_2(d_j)\ge R$. 
\vskip.1cm

\nd It follows from $z_1(d_0)=R$ and definition of $A_1$ that  $d_1<d_0$. Therefore  $0<d_1<d_0\le d_\infty$.
\vskip.1cm 

\nd It remains  to show that   $z_1(d_1)<R$ and $z_2(d_1)=R$.  To do it, let $\{ d_j\} \subseteq  A_1$ such that  $d_j\to d_1$, so that  $z_1(d_j)\to z_1(d_1)\le R$ and $z_2(d_j)\to z_2(d_1)\ge R$.
\vskip.1cm

\nd If $z_1(d_1)=R$ then $u(\cdot,d_1)$ has no zeros in $(0,R)$. By (\ref{T25}), if $d<d_1$ and $d$ is close to $d_1$, $u(\cdot,d)$ has at most one zero in $(0,R)$. If $u(\cdot,d)$ has one zero then we have $d<d_1$ and $d\in A_1$, a contradiction.
\vskip.1cm

\nd On the other hand,  if $u(\cdot,d)$ has  no zero then $d\ge d_0>d_1$ which is again a contradiction. Therefore, $z_1(d_1)<R$. Now, assume by contradiction that $z_2(d_1)>R$. Let $d_j\to d_1$ with $d_j<d_1$. Then, $z_1(d_j)\to z_1(d_1)<R$ and in addition, $z_2(d_j)\to z_2(d_1)>R$, so that, $z_1(d_j)<R$ and $z_2(d_j)>R$ for large $j$, which is impossible. Thus $z_2(d_1)=R$.
\vskip.1cm

\nd By induction, iterating the arguments above, we construct a sequence $\{d_\ell \}_{\ell =1}^{\infty} \subseteq (0,d_\infty]$ such that 
\begin{equation}\label{th12}\begin{split}
0< \cdots <  d_\ell < \cdots < d_1 <  d_\infty,\\
\\ 
z_\ell(d_\ell)<R,\ z_{\ell+1}(d_\ell)=R,
\end{split}\end{equation}

\nd with $d_\ell := \inf A_\ell$, where
$$
A_\ell := \Big \{d\in (0,d_\ell]~|~ \ z_\ell(d)<R,\ z_{\ell+1}(d)\ge R \Big \}.
$$ 
\nd This ends the proof of step 2.
\vskip.2cm

\nd To finish the proof of Theorem \ref{T1}, we use steps 1 and 2 to conclude that for $\lambda\in (0,\Lambda]$ the functions given by Theorem \ref{T2}, namely  $u_\ell(\cdot)=u(\cdot,d_\ell)\in C^1([0,R])$ for $\ell \geq 1$, satisfy

$$r^{\alpha}\phi(|u'_\ell(r)|)u'_\ell(r)\ \mbox{is differentiable},$$

$$-(r^{\alpha}\phi(|u'_\ell(r)|)u'_\ell(r))'=\lambda r^\gamma f(u_\ell(r)),\ 0<r<R,$$

$$u'_\ell(0)=0\ \mbox{and}\ u_\ell(R)=0,$$

\nd that is, $u_l$ is a classical solution of (\ref{P}), $u_\ell$ has precisely $\ell$ zeroes in $(0,R)$ and so 
$$
u_0, u_1, u_2, \cdots,
$$

\nd  is an infinite sequence of solutions of (\ref{P}) as claimed in the statement of Theorem \ref{T1}. \hfill \fbox \hsf

\section{Proof of  Theorem \ref{T2}}

\nd  At first we set
$$
\Phi(t)=\int_0^t s\phi(s)ds,~~  H(t)=t\Phi'(t)-\Phi(t),~~ F(t)=\int_0^t f(s)ds.
$$

\vskip.2cm

\nd The results below will play a crucial role in this paper.

\begin{lem}\label{Prop1} Assume (\ref{alphabeta}) and let $d\in [0,d_\infty]$, $\lambda>0$ and $T>0$. If $u$ is a solution of (\ref{PP}) in $[0,T]$, then

\begin{equation}\label{Prop11} 
H(|u'(r)|)\le \lambda r^{\gamma-\alpha}[F(d)-F(u(r))],\ r\ge 0.
\end{equation}

\begin{equation}\label{Prop12}
 F(u(r))\le F(d)~\mbox{for}~ r\in [0,T] 
\end{equation}
\end{lem}
\vskip.2cm

\begin{lem}\label{T3} Assume that ($\phi_1$)-($\phi_3$), (\ref{alphabeta}) and ($f_1$)-($f_2$) holds. If $f\in C(\mathbb{R})\cap C^1(\mathbb{R}\setminus\{0\})$ and $d \in (0,d_\infty]$,  then problem (\ref{PP}) has a unique solution $u(\cdot,d,\lambda)=u(\cdot,d)=u(\cdot)\in C^1([0,\infty))$.  In addition,

 \begin{equation}\label{T21} \mbox{if}\ d_0  \in (0,d_\infty]\ \mbox{then}\ u(r,d)\to u(r,d_0)\ \mbox{as}\ d\to d_0,\ \mbox{uniformly in}\ [0,T]~ \mbox{for}\ T>0,\end{equation}
  \begin{equation}
\begin{split}
&\mbox{if}\ d_0 \in (0,d_\infty]\ \mbox{then}\ u'(r,d)\to u'(r,d_0)\ \mbox{as}\ d\to d_0, \mbox{uniformly on compact} \\
&\mbox{subsets of}\ (0,\infty],
\end{split}
\label{T22} 
\end{equation}
\end{lem}

\subsection{Proofs of Lemmas \ref{Prop1} and \ref{T3}}\label{Pf lemmata} 

 \nd Integrating the equation in  (\ref{PP})  we get to
 \begin{equation}\label{pr1}
\phi(|u'(r)|)u'(r)=-r^{-\alpha}\int_0^r \lambda t^\gamma f(u(t))dt,\ r>0.
\end{equation}
\nd Setting
\begin{equation}\label{def h}
h(t) := t \phi(t),
\end{equation}
\nd we see that $h$ is invertible with differentiable inverse $h^{-1}$. Then,  

  \begin{equation}\label{pr1h+}
u'(r)=h^{-1}\left(-r^{-\alpha}\int_0^r \lambda t^\gamma f(u(t))dt\right)~\mbox{if}~~ u'(r)>0,
\end{equation}

  \begin{equation}\label{pr1h-}
u'(r)=-h^{-1}\left(-r^{-\alpha}\int_0^r \lambda t^\gamma f(u(t))dt\right)~~\mbox{if}~~  u'(r)<0,
\end{equation}

\nd Once  $f$  is  continuous and $\gamma \geq \alpha$, we conclude from the above equalities that $u\in C^1$.
\vskip.2cm

\nd {\bf  Proof of Lemma \ref{Prop1}}.  From (\ref{pr1h+}) and (\ref{pr1h-}), we infer that   $u\in C^2$ at the points $r>0$ where $u'(r)\neq 0$. Computing derivatives in  (\ref{PP}) and multiplying the resulting equality by $u'(r)$, we are led to
\begin{equation}\label{pr2}
-\alpha r^{\alpha-1}\phi(|u'(r)|)|u'(r)|^2-r^\alpha \frac{d}{dt}h(|u'(r)|)u'(r)u''(r)=\lambda r^{\gamma} f(u(r))u'(r),\ u'(r)\neq 0.
\end{equation}

\nd Consider the functional $E:[0,\infty)\to\mathbb{R}$ defined by $$E(0)=\lambda F(d)\ \mbox{and}\ E(r)=r^{\alpha-\gamma}[H(|u'(r)|)]+\lambda F(u(r)),\ r>0,$$

\nd where $H(t)=t\Phi'(t)-\Phi(t)=\int_0^t h^{\prime}(s)s\ ds$. Note that 

$$
E'(r)=r^{\alpha-\gamma}[H(|u'(r)|)]'+(\alpha-\gamma)r^{\alpha-\gamma-1}H(|u'(r)|)+\lambda f(u(r))u'(r),\ r>0
$$

\nd and $$[H(|u'(t)|)]'=\frac{d}{dt}h(|u'(r)|)u'(r)u''(r),\ u'(r)\neq 0.$$

\nd Therefore from (\ref{pr2}), 
$$
E'(r)=(\alpha-\gamma)r^{\alpha-\gamma-1}H(|u'(r)|)-\alpha r^{\alpha-\gamma-1}\phi(|u'(r)|)|u'(r)|^2\ u'(r)\neq 0.
$$

\nd From Lemma \ref{LA4} in the Appendix, the last inequality combined with hypothesis (\ref{alphabeta}) gives 

 \begin{equation}\label{pr3}E'(r)\le \frac{\gamma_1-1}{\gamma_1}(\alpha-\gamma)r^{\alpha-\gamma-1}\phi(|u'(r)|)|u'(r)|^2-\alpha r^{\alpha-\gamma-1}\phi(|u'(r)|)|u'(r)|^2< 0,\ u'(r)\neq 0.\end{equation}

 \nd Next, we will prove that $E$ is continuous at the origin and therefore, as $E$ is non-decreasing by the previous inequality, it follows  that $E(r)\le E(0)$ for all $r\ge 0$. Note that equation (\ref{pr1}) implies $$\Phi(|u'(r)|)=\Phi\left(h^{-1}\left(\left|r^{-\alpha}\int_0^r \lambda t^{\gamma} f(u(t))dt\right|\right)\right),$$

\nd which in turn gives \begin{equation}\label{pr4}\Phi(|u'(r)|)\le\Phi\left(h^{-1}\left(\frac{\lambda C_{\delta,d} r^{\gamma-\alpha+1}}{\gamma+1 }\right)\right),~ r\in [0,\delta), \end{equation}

\nd where $C_{\delta,d}=\max_{r\in [0,\delta]}|f(u(r))|$. We choose $\delta>0$ small and apply Lemmas \ref{LA1} and \ref{AL} to conclude from (\ref{pr4}) that \begin{equation}\label{pr5} \Phi(|u'(r)|)\le\left(\frac{\lambda C_{\delta,d} }{\gamma+1 }\right)^{\frac{\gamma_2}{\gamma_1-1}}r^{\frac{\gamma_2}{\gamma_1-1}(\gamma-\alpha+1)}, \forall\ r\in [0,\delta).\end{equation}

\nd We apply condition $\Delta _2$ (see  inequality (\ref{delta2}) in the Appendix) in the definition of $E$ to infer that  \begin{equation}\label{pr6}E(r)\le (\gamma_2-1)r^{\alpha-\gamma}\Phi(|u'(r)|)+\lambda F(u(r)),\ r>0.\end{equation}

\nd Thus,  (\ref{pr5})  and  (\ref{pr6}) lead to  
\begin{equation}\label{pr7}
E(r)\le Cr^{(\alpha-\gamma)+\frac{\gamma_2}{\gamma_1-1}(\gamma-\alpha+1)}+\lambda F(u(r)),\ r\in [0,\delta).
\end{equation}

\nd Recalling that  $\gamma_2\ge \gamma_1$, we have that 
$$(
\alpha-\gamma)+\frac{\gamma_2,}{\gamma_1-1}(\gamma-\alpha+1)\ge \frac{\gamma_1-\alpha+\gamma}{\gamma_1-1} > 0.
$$

\nd Hence, from (\ref{pr7}) that $\displaystyle \lim_{r \to 0}E(r)\leq \lambda F(d)$. On the other hand, by Lemma \ref{LA4}, we know that $H(t) \geq 0$ for all $t \geq 0$. Then, by definition of $E$, $E(r) \geq \lambda F(u(r))$ for all $ r >0$. Gathering these information, we conclude that
$$
\lim_{r \to 0}E(r)=\lambda F(d).
$$
Therefore, as (\ref{pr3}) is true,  
$$
E(r)\le E(0)~ \mbox{for}~ r\ge 0,
$$

\nd which is equivalent to the desired inequality namely (\ref{Prop11} )  \hfill \fbox \hsf
\vskip.3cm

\nd  {\bf \large Proof of Lemma \ref{T3} } We will at first study existence and uniqueness of local solutions of  (\ref{PP}).  Let $\epsilon>0$ and consider

\begin{equation}\label{PPP}\tag{$P_{\lambda,d,\epsilon}$}
\left\{\begin{array}{rllr}
-( r^\alpha\phi(|u'(r)|)u'(r))' =&  \lambda r^\gamma f(u(r)),~~ ~~0<r<\epsilon,\\
\\
u(0)=d,&u'(0)=0.
 \end{array}
\right.
\end{equation}
\vskip.3cm

\nd We shall need the following result  whose proof is left to the Appendix.

\begin{lem}\label{LEMA loc sol}
 (\ref{PPP}) has a unique solution $u(\cdot)=u(\cdot,d,\lambda,\epsilon)\in C^2([0,\epsilon))$ for small $\epsilon$.
\end{lem}

\nd {\bf Proof of Uniqueness in Lemma \ref{T3}}  Assume that $u,v$ are two $C^1([0,\infty))$ solutions. Let 
$$
S_0=\{r\ge 0:\ u(t)=v(t),~ 0\le t\le r\}.
$$
\nd We will show  that \begin{equation}\label{th33}S_0\neq\emptyset,\ S_0\ \mbox{is both open and closed in}\ [0,\infty).\end{equation}

\nd Indeed, by Lemma  \ref{LEMA loc sol} above,  $[0,\epsilon)\subset S_0$ for $\epsilon > 0$ small enough,  which shows that $S_0\neq\emptyset$. Moreover, since $u,v$ are $C^1$ functions  we infer  that $S_0$ is  closed. To finish  we shall  prove that $S_0$ is  open.  To achieve that let $\widehat{r}\in S_0$ with $\widehat{r}>0$. We distinguish between two cases.
\vskip.2cm

\nd {\bf{Case 1}}. $u'(\widehat{r})=v'(\widehat{r})=0$
\vskip.3cm
\nd  Assume $u(\widehat{r})=v(\widehat{r})=\widehat{d}$. If  $\widehat{d}=0$ then, up to a translation in the domain, we are within the  settings of Lemma \ref{Prop1}. Therefore,    using  (\ref{Prop11}), observing that   by hypothesis ($f_1$) one has  $F(u(r))\ge 0$ for $r \geq \widehat{r}$, and noticing that   $F(\widehat{d})=0$,  we get
 $$
 H(|u'(r)|)\le \lambda r^{\gamma-\alpha}\left(F(\widehat{d})-F(u(r))\right)\le 0~\mbox{for}~ r \geq  \widehat{r},
 $$
from where it follows that $u(r)=0$ for $r \geq \widehat{r}$, because by Lemma \ref{LA4} in the Appendix 
$$
H(t) \geq 0 \,\,\, \forall t \geq 0 \,\,\, \mbox{and} \,\,\, H(t)=0 \Leftrightarrow t=0.
$$
The same argument works to prove that $v(r)=0$ for $r \geq \widehat{r}$. Consequently, $r\ge \widehat{r}$, $u(r)=v(r)=0$ and then, $S_0=[0,\infty)$ is open. On the other hand, if $\widehat{d}>0$, we define $$\widehat{K}_\rho^\epsilon(\hat{d})=\{u\in C([\widehat{r},\widehat{r}+\epsilon]):\ u(0)=\widehat{d},\ \|u-\widehat{d}\|_\infty\le \rho\},$$

$$\widehat{T}(u(r))=\widehat{d}-\int_{\widehat{r}}^r h^{-1}\left(t^{-\alpha}\int_{\widehat{r}}^t \lambda \tau^\gamma f(u(\tau))d\tau\right)dt,\ \forall\ r\in [0,\epsilon],$$

\nd where $\epsilon,\rho$ are positive  and $\epsilon$ is small. The same proofs of (\ref{th31}) and (\ref{th32}) can be do here, and then the Banach Fixed Point Theorem guarantees a unique fixed point for the operator $\widehat{T}$ when $\epsilon$ is small, therefore, $u(r)=v(r)$ in a small neighbourhood of $\widehat{r}$, which implies that $S_0$ is open.
\vskip.3cm
\nd {\bf{Case 2}}. $u'(\widehat{r})=v'(\widehat{r})\neq 0$.
\vskip.3cm
\nd Note that there is a neighbourhood $V$ of $\widehat{r}$ such that $u'(r),v'(r)\neq 0$ for $r\in V$. So in $V$, we must conclude, as in (\ref{pr3}) (here we use the same notation as in the proof of Lemma \ref{Prop1}) that if $z$ denotes either  $u$ or $ v$ then

\begin{equation*}
\left( r^{\alpha-\gamma}H(|z'(r)|)+\lambda F(z(r))\right)'= -\frac{\alpha+\gamma(\gamma_1-1)}{\gamma_1} r^{\alpha-\gamma-1}\phi(|z'(r)|)z'(r)|^2.
\end{equation*}

\nd Integrating from $\widehat{r}$ to $r$ and subtracting the corresponding equations for $z=u$ and $z=v$, we obtain (remember that $u(\widehat{r})=v(\widehat{r})$ and $u'(\widehat{r})=v'(\widehat{r})$)

\begin{equation}\label{th34}
\begin{split}& r^{\alpha-\gamma}[H(|u'(r)|)-H(|v'(r)|)]+ \lambda \left[F(u(r))-F(v(r))\right]= \\ &-\frac{\alpha+\gamma(\gamma_1-1)}{\gamma_1}
\int_{\widehat{r}}^r t^{\alpha-\gamma-1}  \left[  \phi(|u'(t)|)|u'(t)|^2- \phi(|v'(t)|)|v'(t)|^2\right  ]dt. 
\end{split}\end{equation}

\nd Next  we consider  three  auxliary continuous functions, namely

\begin{equation*}
A_1(t)=\left\{\begin{array}{ccc}
\frac{H(|u'(t)|)-H(|v'(t)|)}{u'(t)-v'(t)}, && u'(t)\neq v'(t)\\
\\
        \phi(|u'(t)|)u'(t),     && u'(t)=v'(t),
 \end{array}
\right.
\end{equation*}

\begin{equation*}
 A_2(t)=\left\{\begin{array}{ccc}
\displaystyle \frac{h(|u'(t)|)|u'(t)|-h(|v'(t)|)|v'(t)|}{u'(t)-v'(t)}, && u'(t)\neq v'(t)\\
\\
        \frac{d}{dt}[h(|u'(t)|)|u'(t)|],     && u'(t)=v'(t),
 \end{array}
\right.
\end{equation*}

\begin{equation*}
B(t)=\left\{\begin{array}{ccc}
\lambda\frac{F(u(t))-F(v(t))}{u(t)-v(t)}, && u(t)\neq v(t)\\
\\
        \lambda f(u(t)),     && u(t)=v(t).
 \end{array}
\right.
\end{equation*}

\nd Let $w(r)=u(r)-v(r)$. From (\ref{th34}), 

\begin{equation}\label{th35}r^{\alpha-\gamma}A_1(r)w'(r)+B(r)w(r)= -\frac{\alpha+\gamma(\gamma_1-1)}{\gamma_1}\int_{\hat{r}}^r t^{\alpha-\gamma-1}A_2(t)w'(t)dt. \end{equation}

\nd Once $u'(\widehat{r})\neq 0$, we have that in a neighbourhood of $\widehat{r}$, the function $1/A_1$ is well defined and continuous,  and so,  equation (\ref{th35}) is the same as

\begin{equation}\label{th36} w'(r)+\frac{B(r)}{A_1(r)}r^{\gamma-\alpha}w(r)=-\frac{\alpha+\gamma(\gamma_1-1)}{\gamma_1}\frac{r^{\gamma-\alpha}}{A_1(r)}\int_{\hat{r}}^r t^{\alpha-\gamma-1}A_2(t)w'(t)dt. \end{equation}

\nd As $h$ is two times differentiable and $u'(\widehat{r})\neq 0$, $A_2$ is continuously differentiable in a neighborhood of $\widehat{r}$, therefore, from (\ref{th36}) and integration by parts we obtain

\begin{equation*}\begin{split}w'(r)+\frac{B(r)}{A_1(r)}r^{\gamma-\alpha}w(r)= \frac{\alpha+\gamma(\gamma_1-1)}{\gamma_1}\frac{r^{\gamma-\alpha}}{A_1(r)}r^{\alpha-\gamma-1}A_2(r)w(r)+ \\  -
 \frac{\alpha+\gamma(\gamma_1-1)}{\gamma_1}\frac{r^{\gamma-\alpha}}{A_1(r)}\int_{\widehat{r}}^r \left[t^{\alpha-\gamma-1}A_2(t)\right]'w(t)dt, \end{split}\end{equation*}

\nd hence
\begin{equation}\label{th37}\begin{split}w'(r)+& \left[\frac{B(r)}{A_1(r)}r^{\gamma-\alpha}-\frac{\alpha+\gamma(\gamma_1-1)}{\gamma_1}\frac{r^{\gamma-\alpha}}{A_1(r)}r^{\alpha-\gamma-1}A_2(r)\right]w(r)= \\ &-
 \frac{\alpha+\gamma(\gamma_1-1)}{\gamma_1}\frac{r^{\gamma-\alpha}}{A_1(r)}\int_{\widehat{r}}^r \left[t^{\alpha-\gamma-1}A_2(t)\right]'w(t)dt. \end{split}\end{equation}

\nd  We introduce the notation
 $$D_1(r)=\frac{B(r)}{A_1(r)}r^{\gamma-\alpha}-\frac{\alpha+\gamma(\gamma_1-1)}{\gamma_1}\frac{r^{\gamma-\alpha}}{A_1(r)}r^{\alpha-\gamma-1}A_2(r),$$

$$D_2(r)=-
 \frac{\alpha+\gamma(\gamma_1-1)}{\gamma_1}\frac{r^{\gamma-\alpha}}{A_1(r)},$$

$$D_3(r)=\left[t^{\alpha-\gamma-1}A_2(t)\right]',$$

\nd which implies from (\ref{th37}) that \begin{equation}\label{th38}w'(r)+D_1(r)w(r)=D_2(r)\int_{\hat{r}}^r D_3(s)w(s)ds.\end{equation}

\nd We integrate equation (\ref{th38}) from $\widehat{r}$ to $r$, which combined with the fact that, $A_1,A_2,1/A_1,A'_2,B$ are bounded functions (remember they are all continuous functions in a neighborhood of $\hat{r}$) to conclude that

\begin{eqnarray*}
 |v(r)| &\leq& \int_{\widehat{r}}^r |D_1(s)||v(s)|ds+\int_{\widehat{r}}^r |D_2(s)|\int_{\widehat{r}}^s |D_3(r)||v(t)|dtds \\
 &\leq& C\int_{\widehat{r}}^r |v(s)|ds,
\end{eqnarray*}
\nd where $C>0$ is a constant. By the Gronwall Inequality, $v=0$ in a neighborhood of $\hat{r}$. Therefore, $S_0$ is open and (\ref{th33}) is proved. 
\vskip.2cm

\nd {\bf Proof of Existence in Lemma \ref{T3}}.~ Let 

\begin{equation*}S_\infty=\{r>0~|~ \mbox{(\ref{PP}) has a solution in}\ [0,r)\}~\mbox{and}~ T_\infty=\sup S_\infty.
\end{equation*}

\nd We will prove that 
 \begin{equation}\label{th311} 
T_\infty=\infty.
\end{equation}

\nd Assume, on the contrary, that $T_\infty <\infty$. First note that $S_\infty$ is a closed set. Indeed, let $r_n\to r$ with $r_n\in S_\infty$. If $r<r_n$ for some $n$ then $r\in S_\infty$ by force, so we can assume that $r_n<r$ and without loss of generality that $r_n<r_{n+1}$. If $u_n$ are the solutions associated with $r_n$, we define $u:[0,r)\to \mathbb{R}$ by $u(x) :=u_n(x)$ of $x\in [0,r_n)$. Once (\ref{th33}) is satisfied, we conclude that $u$ is well defined and it is a solution of (\ref{PP}), which implies that $r\in S_\infty$.
\vskip.2cm

\nd Since $S_\infty$ is closed, we have that $T_\infty\in S_\infty$. Let $u$ be the solution associated with $T_\infty$. We first observe that from (\ref{Prop11}), $|u'|$ is bounded, which implies that $u$ can be continuously extended to $T_\infty$. Moreover, equation (\ref{pr1}) guarantees that $u'(T_\infty)$ is uniquely defined, so there are two cases to consider.:
\vskip.3cm
\nd {\bf{Case 1}}. $u'(T_\infty)=0$.
 \vskip.3cm

\nd If $u(T_\infty)=0$, consider the extension of $u$ namely $\widetilde{u}:[0,\infty)\to \mathbb{R}$ given by $\widetilde{u}(t)=0$ for $t\ge T_\infty$. Then $\widetilde{u}$ is a $C^1$ function and it is also a solution of (\ref{PP}), which is an absurd. Otherwise, if $u(T_\infty)=d^\infty>0$, consider the operator $T$ defined by $$T(u(r))=d^\infty-\int_{T_\infty}^r h^{-1}\left(t^{-\alpha}\int_{T_\infty}^t \lambda \tau^\gamma f(u(\tau))d\tau\right)dt.$$

\nd Following  the same lines as in  the proof of  either (\ref{th31}) or (\ref{th33}) {\bf Case 1}, we have that $T$ has unique fixed point $v:[0,T_\infty+\epsilon]$, which is an absurd due to  the definition of $T_\infty$.
\vskip.3cm

\nd {\bf{Case 2}}. $u'(T_\infty)\neq 0$.
\vskip.1cm

\nd Assume without loss of generality that $u'(T_\infty)>0$. Then, by  (\ref{pr1}),  
$$
u'(r)=h^{-1}\left(r^{-\alpha}\int_0^r\lambda t^\gamma f(u(t))dt\right)
$$
\nd in a neighborhood of $T_\infty$. Hence, $u$ is $C^2$ in a neighborhood of $T_\infty$, which implies by (\ref{PP}) that $$u''(r)=-\left[\frac{d}{dt}h(u'(r))\right]^{-1}\left(\frac{\alpha}{r}\phi(u'(r))u'(r)+\lambda r^{\gamma-\alpha}f(u(r))\right).$$
\vskip.1cm

\nd By the last equation and Peano's Theorem,  $u$ can be extended to  $[0,T_\infty+\delta)$, where $\delta>0$  and thus we reach an absurd, because such extension is also  a solution to (\ref{PP}). This finishes the proof of {\bf Case 2}.  Therefore, (\ref{th311}) is proved and thus Claim 2 is also proved.
\vskip.1cm

 \nd {\bf Proof of (\ref{T21})}. Remember that \begin{equation}\label{th21} r^\alpha\phi(|u'(r)|)u'(r)=-\int_0^r \lambda t\gamma f(u(t))dt.\end{equation} Assume that $d_n\to d_0$ and set $u_n(r)=u(r,d_n)$, $u_0(r)=u(r,d_0)$. Inequality (\ref{pr1}) implies that $|u'_n(r)|$ is bounded for $r\in [0,T]$, therefore, by \'Ascoli-Arz\'ela Theorem,  there is a subesequence, still denoted by $u_n$, such that $u_n\to v$ uniformly in $[0,T]$ for some $v\in C([0,T])$. Now we will prove that $v=u_0$.
\vskip.1cm

\nd First note that by Lebesgue's Theorem

$$\int_0^r \lambda t\gamma f(u_n(t))dt\to \int_0^r \lambda t\gamma f(v(t))dt,$$

\nd and by  (\ref{th21}),

 $$r^\alpha\phi(|u_n'(r)|)u_n'(r)\to -\int_0^r \lambda t\gamma f(v(t))dt,\ r\in [0,T].$$

\nd As a consequence,
  \begin{equation}\label{th22}
|u_n'(r)|\to h^{-1}\left(r^{-\alpha}\left|\int_0^r \lambda t\gamma f(v(t))dt\right|\right),\ r\in [0,T].
\end{equation}

\nd The combination of (\ref{th21}) and (\ref{th22}) implies that $u'_n(r)\to w(r)$ for all $r\in [0,T]$ where $w$ is a continuous function. Hence, applying Lebesgue's Theorem we obtain $$u_n(r)-d_n=\int_0^r u'_n(t)dt\to \int_0^r w(t)dr,\ r\in [0,T],$$

\nd which implies that $w(r)=v'(r)$ and $v'(0)=0$. Once $$\phi(|v'(r)|)v'(r)=-r^{-\alpha}\int _0^r \lambda t^\gamma f(v(t))dt,$$

\nd is satisfied and since $v(0)=d_0$, it follows by the uniqueness of solutions given by theorem \ref{T3} that $v=u_0$, which concludes the proof of (\ref{T21}).
\vskip.2cm

\nd {\bf Proof of (\ref{T22})}. Let $0<a\le r\le b<\infty$ and assume that $d_n\to d_0$. By (\ref{th21}), 
$$
r^{\alpha}|\phi(|u'_n(r)|)u_n'(r)-\phi(|u'_0(r)|)u_0'(r)|\le \int _0^r \lambda t^\gamma |f(u_n(t))-f(u_0(t))|dt.
$$

\nd Since  $(u_n)$ converges uniformly  to $u_0$ in $[a,b]$, we conclude from the previous inequality that $$(\phi(|u'_n(r)|)u_n'(r)-\phi(|u'_0(r)|)u_0'(r))(u'_n(r)-u'_0(r))\to 0,$$

\nd uniformly in $[a,b]$. Now, we combine a generalized form of Simon's inequality, see Lemma \ref{LA3} in the Appendix, with the last convergence to conclude that $u'_n\to u'_0$ uniformly in $[a,b]$. This  finishes the proof of Lemma \ref{T3}. \hfill \fbox \hsf

\subsection  { Proof of Theorem \ref{T2} (Continued)}  

\nd {\bf Proof of (\ref{T23})}. We will start by proving that there is $z_1=z_1(d)>0$ such that $u(z_1)=0$, $u'(z_1)<0$ and \begin{equation}\label{th23} u(r)>0,\ u'(r)<0~ \mbox{for}~0<r<z_1.\end{equation}

\nd Suppose, on the contrary, that $u(r)>0$ for all $r>0$. It follows from (\ref{th21}) and conditions ($f_1$), ($f_2$) that $u'(r)<0$ and

$$-u'(r)\ge h^{-1}\left(\lambda\frac{r^{\gamma-\alpha+1}}{\gamma+1} f(u(r))\right),\ r>0.$$

\nd Note that $u'(r)\to 0$ if $r\to \infty$ because $u(r)>0$.  Hence, the previous inequality implies that $u(r)\to 0$ if $r\to \infty$. Moreover, by  Lemma \ref{LA1} and the previous inequality, we also obtain

$$-u'(r)\ge \max\left\{\left(\frac{\lambda r^{\gamma-\alpha+1} f(u(r))}{(\gamma+1)h(1)}\right)^{\frac{1}{\gamma_1-1}},\left(\frac{\lambda r^{\gamma-\alpha+1} f(u(r))}{(\gamma+1)h(1)}\right)^{\frac{1}{\gamma_2-1}}\right\}, r>0,$$

\nd which implies $$-u'(r)\min\{f(u(r))^{\frac{-1}{\gamma_1-1}},f(u(r))^{\frac{-1}{\gamma_1-2}}\}\ge \min\left\{\left(\frac{\lambda r^{\gamma-\alpha+1}}{(\gamma+1)h(1)}\right)^{\frac{1}{\gamma_1-1}},\left(\frac{\lambda r^{\gamma-\alpha+1}}{(\gamma+1)h(1)}\right)^{\frac{1}{\gamma_2-1}}\right\}$$

\nd for each $r > 0$. We integrate the last inequality from $0$ to $r$ and apply the change of variables $t=u(s)$ to conclude that \begin{equation}\label{th24}\int_{u(r)}^d \min\{f(t)^{\frac{-1}{\gamma_1-1}},f(t)^{\frac{-1}{\gamma_2-1}}\}dt\ge \int _0^r \min\left\{\left(\frac{\lambda s^{\gamma-\alpha+1}}{(\gamma+1)h(1)}\right)^{\frac{1}{\gamma_1-1}},\left(\frac{\lambda s^{\gamma-\alpha+1}}{(\gamma+1)h(1)}\right)^{\frac{1}{\gamma_2-1}}\right\}ds.\end{equation}
\vskip.1cm

\nd Hypothesis (\ref{alphabeta}) implies that the right hand side of (\ref{th24}) converges to infinity as $r\to \infty$. Therefore, (\ref{th24}) yields 
$$
\liminf_{r\to\infty} \int_{u(r)}^d \min\{f(t)^{\frac{-1}{\gamma_1-1}},f(t)^{\frac{-1}{\gamma_2-1}}\}dt=\infty,
$$
\nd which combined with ($f_1$) and the fact that $u(r)\to 0$ if $r\to \infty$, implies a contradiction to  ($f_3$) and thus, (\ref{th23}) is true. To proceed, we will prove that there is $\Lambda>0$ such that \begin{equation}\label{lambda}z_1(d_\infty,\lambda)\ge R\ \mbox{if}\ 0<\lambda\le \Lambda.\end{equation}
\vskip.2cm

\nd Indeed, by (\ref{th21}), 
$$-u'(r)\le h^{-1}\left(\frac{\lambda f(d_\infty) r^{\gamma-\alpha+1}}{\gamma+1}\right)~\mbox{for}~ r\in [0,z_1(d_\infty,\lambda)].
$$

\nd Integrating from $0$ to $r\in [0,d_\infty]$ and making use of Lemma \ref{LA1}, we get that \begin{equation}\label{th25}\begin{split} -u(r)+d_\infty\le \max\Bigg\{&(\gamma_1-1)\left(\frac{\lambda f(d_\infty)}{(\gamma+1)h(1)}\right)^{\frac{1}{\gamma_1-1}}\frac{r^{\frac{\gamma-\alpha+\gamma_1}{\gamma_1-1}}}{\gamma-\alpha+\gamma_1}, \\ &(\gamma_2-1)\left(\frac{\lambda f(d_\infty)}{(\gamma+1)h(1)}\right)^{\frac{1}{\gamma_2-1}}\frac{r^{\frac{\gamma-\alpha+\gamma_2}{\gamma_2-1}}}{\gamma-\alpha+\gamma_2}\Bigg\} . \end{split}\end{equation}
\nd Let $\nu\in (0,1)$.  Choose $r_\infty(\nu)\in (0,z_1(d_\infty,\lambda))$ such that $u(r_\infty(\nu),d_\infty)=\nu d_\infty$. Set $r=r_\infty(\nu)$ in (\ref{th25}) and choose the maximum value in the right hand side of (\ref{th25})  which actually is   $$
(\gamma_1-1)\left(\frac{\lambda f(d_\infty)}{(\gamma+1)h(1)}\right)^{\frac{1}{\gamma_1-1}}\frac{r_\infty(\nu)^{\frac{\gamma-\alpha+\gamma_1}{\gamma_1-1}}}{\gamma-\alpha+\gamma_1}.
$$  
\nd Take  $R>0$ and choose  $\Lambda_\nu>0$ satisfying
\begin{equation}\label{th26}
1-\nu=\left[\left(\frac{\lambda f(d_\infty)}{(\gamma+1)h(1)}\right)^{\frac{1}{\gamma_1-1}}\frac{\gamma_1-1}{\gamma-\alpha+\gamma_1}\right]^{-1}\frac{\Lambda_\nu^{\frac{1}{\gamma_1-1}}R^{\frac{\gamma-\alpha+\gamma_1}{\gamma_1-1}}}{d_\infty}. 
\end{equation}

\nd We infer from (\ref{th25}) and (\ref{th26}) that $$\left[\left(\frac{\lambda f(d_\infty)}{(\gamma+1)h(1)}\right)^{\frac{1}{\gamma_1-1}}\frac{\gamma_1-1}{\gamma-\alpha+\gamma_1}\right]^{-1}\Lambda_\nu^{\frac{1}{\gamma_1-1}}R^{\frac{\gamma-\alpha+\gamma_1}{\gamma_1-1}} \le (\gamma_1-1)\left(\frac{\lambda f(d_\infty)}{(\gamma+1)h(1)}\right)^{\frac{1}{\gamma_1-1}}\frac{r_\infty(\nu)^{\frac{\gamma-\alpha+\gamma_1}{\gamma_1-1}}}{\gamma-\alpha+\gamma_1}, $$

\nd which implies that $$\Lambda_\nu^{\frac{1}{\gamma_1-1}}R^{\frac{\gamma-\alpha+\gamma_1}{\gamma_1-1}}\le \lambda^{\frac{1}{\gamma_1-1}}r_\infty(\nu)^{\frac{\gamma-\alpha+\gamma_1}{\gamma_1-1}}.$$

\nd Hence, 
\begin{equation}\label{th27} 
R\le r_\infty(\nu)\le z_1(d_\infty,\lambda)\ \mbox{if}\  0<\lambda\le \Lambda_\nu.
\end{equation}
\nd To finish the proof of (\ref{lambda}), first note that the maximum of two continuous functions is a continuous function. Therefore, (\ref{th25}) combined with (\ref{th26}) gives 
$$
\Lambda_\nu^{\frac{1}{\eta-1}} \stackrel{\nu \to 0} \longrightarrow \left(\frac{\lambda f(d_\infty)}{(\gamma+1)h(1)}\right)^{\frac{1}{\eta-1}}\frac{\eta-1}{\gamma-\alpha+\eta}\frac{d_\infty}{R^{\frac{\gamma-\alpha+\eta}{\eta-1}}},
$$

\nd where either  $\eta=\gamma_1$ or $\eta=\gamma_2$ depending on whether the maximum in  (\ref{th25}) is assumed at $\gamma_1$ or $\gamma_2$. Note also that $r_\nu(d_\infty)$ is continuous on $\nu$ and $r_\nu(d_\infty)\to z_1(d_\infty,\lambda)$ as $\nu\to 0$, therefore, we conclude from (\ref{th27}) that $$ R\le z_1(d_\infty,\lambda)\ \mbox{if}\  0<\lambda\le \Lambda,$$

\nd where  $$\Lambda :=\frac{\lambda f(d_\infty)}{(\gamma+1)h(1)}\left(\frac{\eta-1}{\gamma-\alpha+\eta}\right)^{\eta-1}\frac{d_\infty^{\eta-1}}{R^{\gamma-\alpha+\eta}}.$$

\nd Now we will show that there is $z_2=z_2(d)>z_1$ such that $u(z_2)=0$, $u'(z_2)>0$ and
\begin{equation}\label{th28}
u(r)<0,\ z_1<r<z_2. 
\end{equation}
\nd In fact, since $u'(z_1)<0$ then, $u'(r)<0$ in a neighborhood of $z_1$. We start by proving that there is $m_1>z_1$ such that $u'(m_1)=0$. Thus, suppose by contradiction that it is not true, i.e. $u'(r)<0$ for all $r>z_1$. We have by (\ref{Prop12}) that

$$\int_0^{u(r)} f(t)dt\le F(d),\ r\ge 0.$$

\nd If there is some sequence $r_n\to \infty$ such that $u(r_n)\to -\infty$ then, by the previous inequality we infer that 
$$\int_{-\infty}^0 f(s)ds=\lim_n \int_{u(r_n)}^0 f(s)ds \ge -F(d),$$
\nd which is impossible, because ($f_1$), ($f_2$) imply that $\int_{-\infty}^0 f(s)ds =-\infty$. Hence, there is $C>0$ such that $$
u(r)\ge -C,\ u'(r)<0,\  \forall r\ge z_1,
$$
\nd and consequently $u(r)\to L$ as $r\to \infty$ for some $L<0$. Now, by (\ref{Prop11}), $$\frac{\Phi(|u'(r)|)}{r^{\gamma-\alpha+1}}\to 0\ \mbox{as}\ r\to \infty,$$

\nd which implies by using the inequality $\Phi(s)\ge c s^2\phi(s)$ that $$\frac{\phi(|u'(r)|)}{r^{\gamma-\alpha+1}}\to 0.$$
\vskip.1cm

\nd On one hand (\ref{th21}) and the previous limits imply that 
$$\frac{1}{r^{\gamma+1}}\int _0^r t^\gamma f(u(t))dt \to 0,$$

\nd and on  the other side, ($f_1$) and L'Hospital rule imply that 
$$
\displaystyle \lim_{r\to \infty}\frac{1}{r^{\gamma+1}}\int _0^r t^\gamma f(u(t))dt =\lim_{r\to \infty} \frac{r^\gamma f(u(r))}{(\gamma+1)r^\gamma}=\frac{f(L)}{\gamma+1}<0,
$$

\nd which is an absurd. Therefore, $u'(m_1)=0$ for some $z_1<m_1$, so that $$u(r)<0\ \mbox{for}\ z_1<r<m_1\ \mbox{and}\ u'(r)<0\ \mbox{for}\ z_1\le r<m_1.$$
\vskip.1cm

\nd Now, taking $r>m_1$,  $r$ close to $m_1$ we have  $$\int_{m_1}^r t^\gamma f(u(t))<0, $$

\nd which implies by (\ref{th21}) that $$u(r)<0,\ u'(r)>0\ \mbox{for all}\ r>m_1,\ r\ \mbox{close to}\ m_1.$$
\vskip.1cm

\nd Assume by contradiction that $u(r)<0$ for  $r>m_1$, so that $u'(r)>0$. Since by ($f_2$) $$-r^\alpha \phi(|u'(r)|)u'(r)=\lambda \int_{m_1}^r t^\gamma f(u(t))dt\le \frac{\lambda f(u(r))}{\gamma+1}(r^{\gamma+1}-m_1^{\gamma+1}),$$

\nd we get by taking $r>\overline{r}=2^{\frac{1}{\gamma+1}}m_1$ above, that $r^{\gamma+1}-m_1^{\gamma+1}>\frac{r^{\gamma+1}}{2}$ and so

$$-r^\alpha\phi(|u'(r)|)u'(r)\le \frac{\lambda f(u(r))}{2(\gamma+1)}r^{\gamma+1},$$

\nd which,  combined with Lemma \ref{LA1} gives

\begin{equation}\label{th29}u'(r)\ge \min\left\{\left(\frac{-\lambda f(u(r))}{2(\gamma+1)}r^{\gamma-\alpha+1}\right)^{\frac{1}{\gamma_1-1}},\left(\frac{-\lambda f(u(r))}{2(\gamma+1)}r^{\gamma-\alpha+1}\right)^{\frac{1}{\gamma_2-1}}\right\},\ r>\overline{r}. \end{equation}

\nd Integrating in (\ref{th29}) from $\overline{r}$ to $r$, we have 
$$
\int _{\overline{r}}^r u'(t)\max\{(-f(u(t)))^{\frac{-1}{\gamma_1-1}},(-f(u(t)))^{\frac{-1}{\gamma_2-1}}\}dt\geq \\
\\
 \int_{\overline{r}}^r \min\left\{\left(\frac{t^{\gamma-\alpha+1}}{2(\gamma+1)}\right)^{\frac{1}{\gamma_1-1}},\left(\frac{t^{\gamma-\alpha+1}}{2(\gamma+1)}\right)^{\frac{1}{\gamma_2-1}}\right\}dt,
$$
\nd for $\ r>\overline{r}$. Making the change of variables $y=u(t)$,
\begin{equation}\label{th210} 
\int_{u(\overline{r})}^{u(r)} \max\{(-f(t))^{\frac{-1}{\gamma_1-1}},(-f(t))^{\frac{-1}{\gamma_2-1}}\}dt\ge \int_{\overline{r}}^r \min\left\{\left(\frac{t^{\gamma-\alpha+1}}{2(\gamma+1)}\right)^{\frac{1}{\gamma_1-1}},\left(\frac{t^{\gamma-\alpha+1}}{2(\gamma+1)}\right)^{\frac{1}{\gamma_2-1}}\right\}dt.
\end{equation}
\nd Once $u(r)<0$ and $u'(r)>0$ for $r>\overline{r}$ it follows that  $u'(r)\to 0$ as $r\to \infty$. Hence, inequality (\ref{th29}) implies that $u(r)\to 0$ as $r\to \infty$. Moreover, the right hand side of (\ref{th210})  converges  to $\infty$ due to  hypothesis (\ref{alphabeta}). Therefore $$\liminf \int_{u(\overline{r})}^{u(r)} \max\{(-f(t))^{\frac{-1}{\gamma_1-1}},(-f(t))^{\frac{-1}{\gamma_2-1}}\}=\infty ,$$

\nd which contradicts ($f_3$), so (\ref{th28}) is proved. Now we will prove that there is $z_3=z_3(d)>z_2$ such that $u(z_3)=0$, $u'(z_3)<0$ and
\begin{equation}\label{th211} u(r)>0\ \mbox{for all}\ r\in (z_2,z_3).\end{equation}
\nd Indeed, since by (\ref{th28}), $u'(z_2)>0$, so that $$u'(r)>0\ \mbox{for all}\ r>z_2,\ r\ \mbox{close to}\ z_2.$$
\vskip.1cm

\nd We claim that there is $m_2>z_2$ such that $u'(m_2)=0$. In fact, othewise, $u'(r)>0$, for all $r>z_2$, which gives that $u(r)>0$  for  $r>z_2$. By (\ref{Prop12}), $$\int_0^{u(r)} f(t)dt\le \int_0 ^d f(t)dt,$$

\nd so that $u(r)\le d$ for  $r\ge z_2$. Hence, there is $L\in (0,d]$ such that $$u(r)\to L\ \mbox{and}\ u(r)\le L,\ r\ge z_2.$$
\vskip.1cm

\nd As in the proof of (\ref{th28}), 

$$\frac{1}{r^{\gamma+1}}\int _0^r t^\gamma f(u(t))dt \to 0,$$

\nd and $$\lim_{r\to \infty}\frac{1}{r^{\gamma+1}}\int _0^r t^\gamma f(u(t))dt =\lim_{r\to \infty} \frac{r^\gamma f(u(r))}{(\gamma+1)r^\gamma}=\frac{f(L)}{\gamma+1}<0,$$

\nd which is an absurd. As a consequence, there is $m_2>z_2$ such that $u'(m_2)=0$ and $u'(r)>0$, $z_2\le r< m_2$, proving the claim. Assume again, by contradiction, that $u(r)>0$ for all $r>m_2$ so that $u'(r)<0$ also for all $r>m_2$. We have (similar to the proof of (\ref{th28}))

$$-r^\alpha \phi(|u'(r)|)u'(r)=\lambda \int_{m_2}^r t^\gamma f(u(t))dt\ge \frac{\lambda f(u(r))}{\gamma+1}(r^{\gamma+1}-m_2^{\gamma+1}).$$

\nd Setting $\overline{r}=2^{\frac{1}{\gamma+1}}m_2$ and taking $r>\overline{r}$,

$$-r^\alpha\phi(|u'(r)|)u'(r)\ge \frac{\lambda f(u(r))}{2(\gamma+1)}r^{\gamma+1},$$

\nd which  combined with (\ref{LA1}) gives
\begin{equation}\label{th212}-u'(r)\ge \min\left\{\left(\frac{\lambda f(u(r))}{2(\gamma+1)}r^{\gamma-\alpha+1}\right)^{\frac{1}{\gamma_1-1}},\left(\frac{\lambda f(u(r))}{2(\gamma+1)}r^{\gamma-\alpha+1}\right)^{\frac{1}{\gamma_2-1}}\right\},\ r>\overline{r}. 
\end{equation}
\nd Integrating (\ref{th212}) from $\overline{r}$ to $r$ and making the change of variables $u(t)=s$, we get $$\int_{u(\overline{r})}^{u(r)} -\max\{f(t)^{\frac{-1}{\gamma_1-1}},f(t)^{\frac{-1}{\gamma_2-1}}\}dt\ge \int_{\overline{r}}^r \min\left\{\left(\frac{t^{\gamma-\alpha+1}}{2(\gamma+1)}\right)^{\frac{1}{\gamma_1-1}},\left(\frac{t^{\gamma-\alpha+1}}{2(\gamma+1)}\right)^{\frac{1}{\gamma_2-1}}\right\}dt.$$
\nd Taking $\liminf$ in both sides, we arrive at a contradiction with ($f_3$) and so  (\ref{th211}) is true. To finish the proof of (\ref{T22}) we argue as in (\ref{th28}) and (\ref{th211}) to get zeroes $z_4,z_5$ and inductively, a sequence with the properties asserted in (\ref{T23}).
\vskip.3cm

\nd {\bf Proof of (\ref{T24})}. We start by proving that $z_1(d)\to 0$ when $d\to 0$. By (\ref{th21}) and (\ref{th23}) we obtain $$-u'(r)=h^{-1}\left(r^{-\alpha}\int_0^r \lambda t^\gamma f(u(t))dt\right),\ r\in [0,z_1].$$

\nd Now we apply ($f_2$) and Lemma \ref{LA1} to conclude that $$-u'(r)\ge \min\left\{\left(\lambda\frac{r^{\gamma-\alpha+1}f(u(r))}{\gamma+1}\right)^{\frac{1}{\gamma_1-1}},\left(\lambda\frac{r^{\gamma-\alpha+1}f(u(r))}{\gamma+1}\right)^{\frac{1}{\gamma_2-1}}\right\},\ r\in [0,z_1],$$

\nd which implies that 
$$-u'(r)\max\left\{f(u(r))^{\frac{-1}{\gamma_1-1}},f(u(r))^{\frac{-1}{\gamma_2-1}}\right\}\ge \min\left\{\left(\lambda\frac{r^{\gamma-\alpha+1}}{\gamma+1}\right)^{\frac{1}{\gamma_1-1}},\left(\lambda\frac{r^{\gamma-\alpha+1}}{\gamma+1}\right)^{\frac{1}{\gamma_2-1}}\right\},\ r\in [0,z_1]. 
$$
\nd Integrating from $0$ to $r$ and making the change of variables $y=u(t)$ we get to 
$$\int_{u(r)}^d \max\left\{f(t)^{\frac{-1}{\gamma_1-1}},f(t)^{\frac{-1}{\gamma_2-1}}\right\}dt\ge \int_0^r \min\left\{\left(\lambda\frac{t^{\gamma-\alpha+1}}{\gamma+1}\right)^{\frac{1}{\gamma_1-1}},\left(\lambda\frac{t^{\gamma-\alpha+1}}{\gamma+1}\right)^{\frac{1}{\gamma_2-1}}\right\}dt.
$$
\nd Taking $r=z_1(d)$ in the previous inequality and making use of (\ref{alphabeta}) and ($f_3$), we  conclude that $z_1(d)\to 0$ as $d\to 0$. Now, letting $\ell \ge 1$, we assume that $u(r)>0$ in $(z_{\ell}(d),z_{\ell +1}(d))$, so that by the notations of (\ref{th28}) and (\ref{th211}) we have $u'(r)>0$ in $(z_{\ell}(d),m_{\ell}(d))$ and $u'(r)<0$ in $(m_l(d),z_{\ell+1}(d))$ (the case $u(r)<0$ in $(z_{\ell}(d),z_{\ell+1}(d))$ is handled similarly). Now, using ($f_2$) in (\ref{th21}), taking $m_{\ell}(d)\le r \le z_{\ell+1}(d)$ and then applying lemma \ref{LA1}, we obtain successively

$$r^\alpha h(-u'(r))\ge \lambda f(u(r))\frac{r^{\gamma+1}-m_{\ell}(d)^{\gamma+1}}{\gamma+1},$$

$$
-u'(r)\max\{f(u(r))^{\frac{-1}{\gamma_1-1}},f(u(r))^{\frac{-1}{\gamma_2-1}}\}\ge \min\left\{\left(\lambda\frac{r^{\gamma+1}-m_{\ell}(d)^{\gamma+1}}{(\gamma+1)r^\alpha}\right)^{\frac{1}{\gamma_1-1}},\left(\lambda\frac{r^{\gamma+1}-m_{\ell}(d)^{\gamma+1}}{(\gamma+1)r^\alpha}\right)^{\frac{1}{\gamma_2-1}}\right\}.
$$
\nd Note that $r^{\gamma-\alpha}\ge m_{\ell}(d)^{\gamma-\alpha}$ since $\gamma\ge\alpha$, therefore 
$$
r^{\gamma-\alpha+1}-r^{-\alpha}m_{\ell}(d)^{\gamma+1}\ge m_{\ell}(d)^{\gamma-\alpha}(r-m_{\ell}(d)),
$$
\nd which gives 
\begin{equation}
\begin{split}
& -u'(r)\max\{f(u(r))^{\frac{-1}{\gamma_1-1}},f(u(r))^{\frac{-1}{\gamma_2-1}}\} \ge \\ & \min\left\{\left[\frac{\lambda m_{\ell}(d)^{\gamma-\alpha}}{(\gamma+1)}(r-m_{\ell}(d))\right]^{\frac{1}{\gamma_1-1}},\left[\frac{\lambda m_{\ell}(d)^{\gamma-\alpha}}{(\gamma+1)}(r-m_{\ell}(d))\right]^{\frac{1}{\gamma_2-1}}\right\}.
\end{split}
\end{equation}
\nd Integrating from $m_{\ell}(d)$ to $z_{\ell+1}(d)$, making the change of variables $y=u(t)$, we find that \begin{equation}\label{th213} \begin{split} &\int_0^{u(m_{\ell}(d))} \max\{f(t)^{\frac{-1}{\gamma_1-1}},f(t)^{\frac{-1}{\gamma_2-1}}\}dt \ge \\ & \int _{m_{\ell}(d)}^{z_{\ell+1}(d)}\min\left\{\left[\frac{\lambda m_{\ell}(d)^{\gamma-\alpha}}{(\gamma+1)}(r-m_{\ell}(d))\right]^{\frac{1}{\gamma_1-1}},\left[\frac{\lambda m_{\ell}(d)^{\gamma-\alpha}}{(\gamma+1)}(r-m_{\ell}(d))\right]^{\frac{1}{\gamma_2-1}}\right\}dt.
\end{split}
\end{equation}

\nd Assume now $z_{\ell}(d)<r<m_{\ell}(d)$. Then by a similar argument, this time, integrating from $z_{\ell}(d)$ to $m_{\ell}(d)$ we deduce that

\begin{equation}\label{th214} \begin{split} &\int_0^{u(m_{\ell}(d))} \max\{f(t)^{\frac{-1}{\gamma_1-1}},f(t)^{\frac{-1}{\gamma_2-1}}\}dt \ge \\ & \int _{z_{\ell}(d)}^{m_{\ell}(d)}\min\left\{\left[\frac{\lambda m_{\ell}(d)^{\gamma-\alpha}}{(\gamma+1)}(m_{\ell}(d)-r)\right]^{\frac{1}{\gamma_1-1}},\left[\frac{\lambda m_{\ell}(d)^{\gamma-\alpha}}{(\gamma+1)}(m_{\ell}(d)-r)\right]^{\frac{1}{\gamma_2-1}}\right\}dt.\end{split}\end{equation}

\nd Now, since $u(m_{\ell}(d))\le d$ we have by ($f_3$) that the left hand side of (\ref{th213}) and (\ref{th214}) converge to zero, and therefore, $\displaystyle \lim_{d\to 0}z_{\ell}(d)=\lim_{d\to 0}z_{\ell+1}(d)$ for each $\ell\geq 1$.  Once $z_1(d)\to 0$ as $d\to 0$, we see that $z_{\ell}(d)\to 0$ as $d\to 0$.
\vskip.3cm

\nd We pass to the proof that $z_{\ell}(d)\to z_{\ell}(d_0)$ if $d\to d_0$. Let us first show that $z_1(d)\to z_1(d_0)$ as $d\to d_0$. Indeed, let $d_n\to d_0$, $u_n(\cdot)=u(\cdot,d_n)$ and $u_0(\cdot)=u(\cdot,d_0)$ so that we have from (\ref{T21}) that $u_n\to u$ uniformly in compact subsets of $(0,\infty)$. For each $\epsilon>0$ small we find

$$u_0(r)>0,\ 0\le r\le z_1(d_0)-\epsilon\ \mbox{and}\ u_0(z_1(d_0)+\epsilon)<0,$$

\nd so that $$u_n(r)>0,\ 0\le r\le z_1(d_0)-\epsilon\ \mbox{and}\ u_n(z_1(d_0)+\epsilon)<0,$$
\vskip.1cm

\nd for sufficiently large $n$. As a consequence, $z_1(d_0)-\epsilon<z_1(d_n)<z_1(d_0)+\epsilon$, showing that $z_1(d_n)\to z_1(d_0)$. Now, assume by induction that $z_{\ell}(d_n)\to z_{\ell}(d_0)$ for some $\ell>1$. We will show that $z_{\ell+1}(d_n)\to z_{\ell+1}(d_0)$. For that matter, we assume $u_0(t)<0$ for $z_{\ell}(d_0)<t<z_{\ell+1}(d_0)$ (the other case is handled similarly). Taking $\epsilon>0$ small, we find that $u_n(t)<0$ for $z_{\ell}(d_0)+\epsilon\le t\le z_{\ell+1}(d_0)-\epsilon$ and $u_n(z_{\ell+1}(d_0)+\epsilon)>0$, showing that $z_{\ell+1}(d_0)-\epsilon<z_{\ell+1}(d_n)<z_{\ell+1}(d_0)+\epsilon$. Consequently, $z_{\ell+1}(d_n)\to z_{\ell+1}(d_0)$ as $d\to d_0$, which finishes the proof of (\ref{T24}).
\vskip.1cm
\nd {\bf Proof of (\ref{T25})}.
\vskip.1cm

\nd Let $d\in (0, d_0)$. It suffices to show that $z_{\ell+2}(d)>R$ whenever $d$ is close enough to $d_0$. We assume that $u(r,d_0)<0$ for $r\in (z_{\ell}(d_0),z_{\ell+1}(d_0))$ (the other case is handled similarly).
\vskip.3cm
\nd Notice that as $z_{\ell}(d_0)$ is increasing an there is only $\ell$ zeroes in $(0,R)$, we must show that $z_{\ell+1}(d_0)\ge R$ and $z_{\ell+2}(d_0)>R$. However, as $z_{\ell+2}(d)\to z_{\ell+2}(d_0)$ for $d\to d_0$, we have $z_{\ell+2}(d)>R$ whenever $d$ is close enough to $d_0$. This completes the proof of Theorem \ref{T2}. \hfill \fbox \hsf

\section{Appendix}

\begin{rmk}
\nd ({\bf On the radially symetric form of  $(\Phi)$})~ Let $u$ be a weak solution of $(\Phi)$, radially symmetric  in the sense that   $u(x)=u(|x|)= u(r)$. Let $r \in (0,R)$ and pick $\epsilon > 0$ small such that  $0 < r < r + \epsilon < R$.
\vskip.2cm

\nd Consider the radially symmetric cut-off function $v_{r,\epsilon}(x) = v_{r,\epsilon}(r) $, where
$$
v_{r,\epsilon}(t) := \left\{ \begin{array}{l}
1~~ \mbox{if}~~ 0 \leq t \leq r,\\
linear~~ \mbox{if}~~ r \leq t \leq r + \epsilon,\\
0~~ \mbox{if}~~ r+\epsilon \leq t \leq R.
\end{array} \right.
$$
\nd and notice that  $v_{r,\epsilon} \in W_0^{1, \Phi}({B})| \cap Lip( {\overline{B}})$. By replacing  $v$  with $v_{r,\epsilon}$ in $(\ref{radialphi})$, we get to
 $$
 \frac{-1}{\epsilon}\int_{B(0,r+\epsilon)\setminus B(0,r)}\phi (|u'(|x|)|)u'(|x|) dx=\lambda\int_{B(0,r+\epsilon)}f(u(|x|))v_{r,\epsilon}(|x|) dx.
 $$

\nd Making the change of variables $x=r\omega$ with $r>0$ and $\omega \in \partial B(0,1)$ and letting  $\epsilon\to 0$ we infer that
$$
\phi(|u'(r)|)u'(r)r^{N-1} = \lambda\int_0^r f(u(r))r^{N-1}dr,
$$

\nd which gives
$$
(r^{N-1}\phi(|u'(r)|)u'(r))'=\lambda r^{N-1}f(u(r)).
$$
\nd So the radially symetric form of   $(\Phi)$    is
$$
\left\{\begin{array}{rllr}
\displaystyle - (r^{N-1}\phi(|u'(r)|)u'(r))' &=& \lambda r^{N-1}f(u(r),~ 0 <  r < R\\
\\
 u'(0)  =  u(R) = 0.
\end{array}
\right.
$$
\end{rmk}

\begin{lem}\label{LA1} Assume that $\phi$ satisfies ($\phi_1$)-($\phi_3$). Then $$h(1)\min\{h^{-1}(s)^{\gamma_1-1},h^{-1}(s)^{\gamma_2-1}\}\le s\le h(1) \max\{h^{-1}(s)^{\gamma_1-1},h^{-1}(s)^{\gamma_2-1}\},\  s>0.$$\end{lem}

\nd \proof. Condition ($\phi_3$) implies that $$(\gamma_1-1)\frac{d}{dt}\ln{t}\le \frac{d}{dt}\ln {h(t)}\le (\gamma_2-1)\frac{d}{dt}\ln{t},\ \forall\ t>0.$$

\nd Let  $t\le 1$. Integrating the previous inequality from $t$ to $1$, we  get $$h(1) t^{\gamma_1-1}\le h(t)\le h(1) t^{\gamma_2-1},\  t\le 1.$$

\nd Let  $t\ge 1$.  Iintegrating the previous inequality from $1$ to $t$, we  get $$h(1) t^{\gamma_2-1}\le h(t)\le h(1) t^{\gamma_1-1},\ \forall\ t\ge 1.$$

\nd Therefore

$$h(1) \min\{t^{\gamma_1-1},t^{\gamma_2-1}\}\le h(t)\le h(1) \max\{t^{\gamma_1-1},t^{\gamma_2-1}\},\ \forall\ t>0.$$

\nd Letting $t=h^{-1}(s)$, the lemma is proved. \hfill \fbox \hsf

\begin{lem}\label{AL} 
Assume $\phi$ satisfies ($\phi_1$)-($\phi_3$). Then $$\Phi(1)\min\{t^{\gamma_1},t^{\gamma_2}\}\le \Phi(t) \le \Phi(1) \max\{t^{\gamma_1},t^{\gamma_2}\},\  t>0.$$
\end{lem}

\nd \proof. From ($\phi_3$),  
$$
\gamma_1 t\phi(t) \le th'(t)+t\phi(t)\le \gamma_2 t\phi(t),\ \forall t>0,
$$

\nd which implies, after  integration from $0$ to $t$ that, 
 \begin{equation}\label{delta2}
\gamma_1\le \frac{t\Phi'(t)}{\Phi(t)}\le \gamma_2,\  t>0.
\end{equation}

\nd The previous inequality is called condition $\Delta_2$. To finish the proof, we proceed as in the proof of lemma \ref{LA1} to conclude the desired inequality. \hfill \fbox \hsf

\begin{lem}\label{ALL} Assume that $\phi$ satisfies ($\phi_1$)-($\phi_3$). Then $$[h^{-1}]'(t)\le \frac{t^{\frac{-\gamma_2+2}{\gamma_2-1}}}{h(1)^{\gamma_2}(\gamma_1-1)},\  t\le 1.$$\end{lem}

\nd \proof. Remember that \begin{equation}\label{APE1}[h^{-1}]'(t)=\frac{1}{h'(h^{-1}(t))},\  t>0.\end{equation}

\nd From the proofs of Lemmas \ref{LA1} and \ref{AL}, 

\begin{equation}\label{APE2}h(1)(\gamma_1-1)\min\{t^{\gamma_1-2},t^{\gamma_2-2}\}\le h'(t)\le h(1)(\gamma_2-1)\max\{t^{\gamma_1-2},t^{\gamma_2-2}\}~\mbox{for}~  t>0.\end{equation}

\nd Gathering (\ref{APE1}) and (\ref{APE2}), we see that 
$$
[h^{-1}]'(t)\le \frac{[h^{-1}(t)]^{-\gamma_2+2}}{h(1)(\gamma_1-1)},\  t\le 1.
$$

\nd Now we use Lemma \ref{LA1} to obtain $$[h^{-1}]'(t)\le \frac{t^{\frac{-\gamma_2+2}{\gamma_2-1}}}{h(1)^{\gamma_2}(\gamma_1-1)},\  t\le 1.$$
\hfill \fbox \hsf

\begin{lem}Assume that $\phi:(0,\infty)\to (0,\infty)$ is a differentiable function satisfying ($\phi_3$). Then, there is a positive constant $\Gamma_1$ such that \begin{equation}\label{LA2} \sum_{i,j=1}^N\frac{\partial a_j}{\partial\eta_i}(\eta)\xi_i\xi_j\geq\Gamma_1\phi(|\eta|)|\xi|^2, \end{equation}
where $a_j(\eta)=\phi(|\eta|)\eta_j$, $\eta\in \mathbb{R}^N\setminus\{0\}$ and $\xi\in \mathbb{R}^N$.\end{lem}

\nd \proof.  Indeed, by ($\phi_3$), 
\begin{equation}\label{1.2}
(\gamma_1-2)\phi(t)\leq t\phi'(t)\leq (\gamma_2-2)\phi(t).
\end{equation} 
Suppose first that $\gamma_1<2$. Note that \begin{equation}\label{1.3}\sum_{i,j=1}^N\frac{\partial a_j}{\partial\eta_i}(\eta)\xi_i\xi_j=\phi(|\eta|)|\xi|^2+\frac{\phi'(|\eta|)|\langle\eta,\xi\rangle|^2}{|\eta|}\end{equation}

\nd If $\phi'(|\eta|)<0$, then $\phi'(|\eta|)|\langle\eta,\xi\rangle|^2\geq\phi'(|\eta|)|\eta|^2|\xi|^2$. From $(\ref{1.2})$ and $(\ref{1.3})$, 
$$
\sum_{i,j=1}^N\frac{\partial a_j}{\partial\eta_i}(\eta)\xi_i\xi_j\geq (\gamma_1-1)\phi(|\eta|)|\xi|^2.
$$

\nd If $\phi'(|\eta|)\geq 0$, then take $\Gamma_1=1$.
\vskip.3cm

\nd If $\gamma_1\geq 2$, then (\ref{1.2}) is satisfied with $\Gamma_1=1$, as can readily be seen from $(\ref{1.3})$ and noting that $\phi'(t)\geq 0$ in this case.  \hfill \fbox \hsf
\vskip.3cm
\nd We now prove a Simon type inequality.

\begin{lem}\label{LA3} Assume that $\phi:(0,\infty)\to (0,\infty)$ is a differentiable function satisfying ($\phi_1$)-($\phi_3$). Then \begin{equation}\label{simoninequality} \langle \phi(|\eta|)\eta-\phi(|\eta'|)\eta',\eta-\eta'\rangle\ge \min\{4,4\Gamma_1\}\frac{|\eta-\eta'|}{1+|\eta|+|\eta'|}\Phi\left(\frac{|\eta-\eta'|}{4}\right),\end{equation} where $\Gamma_1$ was given in lem.a \ref{LA2}, $\eta,\eta'\in \mathbb{R}^N$ and $\langle\cdot,\cdot\rangle$ denotes inner product.\end{lem}
\vskip.1cm

\nd \proof. If $\eta,\eta'=0$ then (\ref{simoninequality}) is obviously satisfied. If only one of them is $0$, let's say, $\eta'=0$, then $$\phi(|\eta|)|\eta|^2\ge \Phi(|\eta|)\ge 4\Phi\left(\frac{|\eta|}{4}\right),$$

\nd where in the last inequalities we have used the properties of an N-function (note that an N-function is convex). So (\ref{simoninequality}) is satisfied. If $\eta,\eta'\neq 0$, assume without loss of generality that $|\eta|\le |\eta'|$. Then, an application of Cauchy-Schwartz inequality implies that 
$$\frac{|\eta-\eta'|}{4}\le |t\eta+(1-t)\eta'|\le 1+|\eta|+|\eta'|,\  t\in [0,1/4].
$$

\nd We conclude from the last inequality, (\ref{1.3}) and the properties of an  N-function that
\begin{eqnarray*}
  \langle \phi(|\eta|)\eta-\phi(|\eta'|)\eta',\eta-\eta'\rangle &=& \sum_{i=1}^N \int_0^1 \frac{d}{dt}[a_j(t\eta+(1-t)\eta')](\eta_j-\eta'_j)dt \\
   &=& \int_0^1\sum _{i,j=1}^N \frac{\partial a_j}{\partial\eta_i}[t\eta+(1-t)\eta'](\eta_i-\eta_i')(\eta_j-\eta'_j)dt \\
   &\ge& \Gamma_1\int_0^1 \phi(|t\eta+(1-t)\eta'|)|\eta-\eta'|^2dt \\
   &\ge& \Gamma_1\int_0^{1/4} \phi(|t\eta+(1-t)\eta'|)|\eta-\eta'|^2dt \\
   &=& \Gamma_1\int_0^{1/4} \phi(|t\eta+(1-t)\eta'|)|\eta-\eta'|^2\frac{|t\eta+(1-t)\eta'|}{|t\eta+(1-t)\eta'|}dt \\
   &\ge& 4\Gamma_1\frac{|\eta-\eta'|}{1+|\eta|+|\eta'|}\phi\left(\frac{|\eta-\eta'|}{4}\right)\left(\frac{|\eta-\eta'|}{4}\right)^2 \\
   &\ge& 4\Gamma_1 \frac{|\eta-\eta'|}{1+|\eta|+|\eta'|} \Phi\left(\frac{|\eta-\eta'|}{4}\right).
\end{eqnarray*}\hfill \fbox \hsf

\begin{lem}\label{LA4}  Assume $\phi$ satisfies ($\phi_1$)-($\phi_3$). Then, the function $H(t)=t\Phi'(t)-\Phi(t)$ is  strictly increasing  and satisfies 
$$
(\gamma_1-1)\Phi(t)\le H(t)\le (\gamma_2-1)\Phi(t),\ t\ge 0,
$$  

$$
\frac{\gamma_1-1}{\gamma_1}t\Phi'(t)\le H(t)\le \frac{\gamma_2-1}{\gamma_2}t\Phi'(t),\ t\ge 0.
$$\end{lem}

\nd \proof. Indeed, as ($\phi_3$) is satisfied, we have that $$t\Phi'(t)-r\Phi'(r)>(t-r)\Phi'(t)>\int_ r^t \tau\phi(\tau) d\tau, t>r\ge 0,$$

\nd which implies that $H$ is strictly incresing. On the other hand, condition (\ref{delta2}) implies the desired inequalities.

\hfill \fbox \hsf

\nd {\bf Proof of Lemma \ref{LEMA loc sol}} Indeed, take $\rho\in (0,d)$ and set  
$$
K_\rho^\epsilon(d)=\{u\in C([0,\epsilon])~|~ u(0)=d,\ \|u-d\|_\infty\le \rho\}.
$$
\nd Take $\epsilon > 0$ small.   If $u\in K_\rho^\epsilon(d)$, then by continuity, $u(r)>0$, $r\in [0,\epsilon]$. Hence, for small $\epsilon$,  a solution of (\ref{PPP}) satisfies $u'(r)\le 0$ for $r\in [0,\epsilon]$ (this was showed in the proof of proposition (\ref{Prop1})) and $$u(r)=d-\int_0^r h^{-1}\left(t^{-\alpha}\int_0^t \lambda \tau^\gamma f(u(\tau))d\tau\right)dt,\ \forall\ r\in [0,\epsilon].$$

\nd We infer  that the solutions of (\ref{PPP}), for small $\epsilon$, are fixed points of the operator $$T(u(r))=d-\int_0^r h^{-1}\left(t^{-\alpha}\int_0^t \lambda \tau^\gamma f(u(\tau))d\tau\right)dt,\ \forall\ r\in [0,\epsilon].$$

\nd Now we will verify that there exist $\epsilon,\rho>0$ and $k\in (0,1)$ such that 
\begin{equation}\label{th31} 
T\left(K_\rho^\epsilon(d)\right)\subset K_\rho^\epsilon(d),
\end{equation}
and
\begin{equation}\label{th32} 
\|Tu-Tv\|_\infty\le k\|u-v\|_\infty.
\end{equation}

\nd Therefore,  by the Banach Fixed Point Theorem, $T$ has a unique fixed point, which in turn will be a $C^2([0,\epsilon])$ solution of (\ref{PPP}). With respect to (\ref{th31}), let $\rho\in (0,d/2]$, which implies that $u(r)\in [d/2,2d]$ for $u\in K_\rho^\epsilon(d)$. Therefore, for $u\in K_\rho^\epsilon(d)$ we have that $$h^{-1}\left(r^{-\alpha}\int_0^s \lambda t^{\gamma} f(u(t))dt \right)\leq h^{-1}\left(\frac{\lambda \|f\|_{\infty,d}s^{\gamma-\alpha+1}}{\gamma+1}\right),\ s\in [0,\epsilon],$$

\nd where $\|f\|_{\infty,d}=\displaystyle \max_{s\in [d/2,2d]}f(s)$. For small $\epsilon$, we can apply lemma \ref{LA1} in the Appendix to conclude from the previous inequality that
\begin{eqnarray*}
  |T(u(r))-T(u(0))| &=& \int_0^r h^{-1}\left(r^{-\alpha}\int_0^s \lambda t^{\gamma}f(u(t))dt\right)ds \\
   &\le& \int_0^rh^{-1}\left(\frac{\lambda \|f\|_{\infty,d}s^{\gamma-\alpha+1}}{\gamma+1}\right) ds  \\
   &\le& \int_0^r h(1) \left( \frac{\lambda \|f\|_{\infty,d}s^{\gamma-\alpha+1}}{\gamma+1}\right)^{\frac{1}{\gamma_1-1}}ds \\
   &=& h(1)\left(\frac{\lambda \|f\|_{\infty,d}}{\gamma+1}\right)^{\frac{1}{\gamma_1-1}}r^{\frac{\gamma-\alpha+\gamma_1}{\gamma_1-1}},\  r\in [0,\epsilon].
\end{eqnarray*}
\nd As $\gamma\ge \alpha$, we obtain from the last inequality that there is $\epsilon>0$ such that $Tu\in C([0,\epsilon])$ and $|T(u(r))-d|\le\rho$ for $r\in [0,\epsilon]$, which finishes the proof  of (\ref{th31}). Now we pass to the proof of (\ref{th32}). We first prove it by assuming that $f\in C(\mathbb{R})\cap C^1(\mathbb{R}\setminus\{0\})$. Fix $\rho$ as in (\ref{th31}) and take $u,v\in K_\rho^\epsilon(d)$. By the Mean Value Theorem, there is $h\in (0,1)$ such that
\begin{equation*}
  \begin{split}T(v(r))-T(u(r))=\int_0^r \left[h^{-1}\left(s^{-\alpha}\int_0^s \lambda t^{\gamma} f(u(t))dt\right)-h^{-1}\left(s^{-\alpha}\int_0^s \lambda t^{\gamma} f(v(t))dt\right)\right]ds= \\ \int_0^r \left[(h^{-1})'\left(s^{-\alpha}\int_0^s \lambda t^{\gamma} f(hu(t)+(1-h)v(t))dt\right)\left(s^{-\alpha}\int_0^s \lambda t^{\gamma} f'(hu(t)+(1-h)v(t))(u(t)-v(t))dt\right)\right]ds.
  \end{split}
\end{equation*}

\nd Choose $\epsilon$ small in such a way that the number $s^{-\alpha}\int_0^s \lambda t^{\gamma} f(hu(t)+(1-h)v(t))dt$ for $s\in [0,\epsilon]$ is small. Therefore, Lemma \ref{ALL} and the last equality implies that for $1<\gamma_2\le 2$ (note that in this case, the function $t\mapsto t^{\frac{-\gamma_2+2}{\gamma_2-1}}$ is decreasing)

\begin{eqnarray*}
  &|T(v(r))-T(u(r))|\le&  \\
  &\int_0^r \left[c\left(s^{-\alpha}\int_0^s \lambda t^{\gamma} |f(hu(t)+(1-h)v(t))|dt\right)^{\frac{-\gamma_2+2}{\gamma_2-1}}\left(s^{-\alpha}\int_0^s \lambda t^{\gamma} |f'(hu(t)+(1-h)v(t))||u(t)-v(t)|dt\right)\right]ds\le&  \\
   &\displaystyle\int_0^r \left[c\left(s^{-\alpha}\int_0^s \lambda \|f\|_{\infty,d'} t^{\gamma}dt\right)^{\frac{-\gamma_2+2}{\gamma_2-1}}\left(s^{-\alpha}\int_0^s \lambda \|f'\|_{\infty,d} t^{\gamma} dt\right)\|u-v\|_\infty\right]ds=& \\
   &\displaystyle\int_0^r \left[c\left(\lambda \|f\|_{\infty,d'} \frac{s^{-\alpha+\gamma+1}}{\gamma+1}\right)^{\frac{-\gamma_2+2}{\gamma_2-1}}\left(\lambda \|f'\|_{\infty,d}\frac{s^{-\alpha+\gamma+1}}{\gamma+1} dt\right)\|u-v\|_\infty \right]ds=& \\
   &\displaystyle c\left(\frac{\lambda\|f\|_{\infty,d'}}{\gamma+1}\right)^{\frac{-\gamma_2+2}{\gamma_2-1}}\frac{\lambda\|f'\|_{\infty,d}}{\gamma+1} r^{\frac{-\alpha+\gamma+\gamma_2}{\gamma_2-1}}\|u-v\|_\infty,&
\end{eqnarray*}

  \nd where $\|f\|_{\infty,d'}=\min_{s\in [d/2,2d]}|f(s)|$ and $\|f'\|_{\infty,d}=\max_{s\in [d/2,2d]}|f'(s)|$. If on the other hand, we have that $\gamma_2\ge 2$, i.e., $t\mapsto t^{\frac{-\gamma_2+2}{\gamma_2-1}}$ is increasing then, we must conclude that

  $$|T(v(r))-T(u(r))|\le c\left(\frac{\lambda\|f\|_{\infty,d}}{\gamma+1}\right)^{\frac{-\gamma_2+2}{\gamma_2-1}}\frac{\lambda\|f'\|_{\infty,d}}{\gamma+1} r^{\frac{-\alpha+\gamma+\gamma_2}{\gamma_2-1}}\|u-v\|_\infty,$$

  \nd where $\|f\|_{\infty,d}=\max_{\in [d/2,d]}f(s)$. In both cases, hypothesis (\ref{alphabeta})  implies the existence of $\epsilon$ such that (\ref{th32}) is true in the case $f\in C(\mathbb{R})\cap C^1(\mathbb{R}\setminus\{0\})$.

\hfill \fbox \hsf

\begin{flushright}

\scriptsize{ \bf Claudianor O. Alves   }\\

  \scriptsize{Universidade Federal de Campina Grande\\
   Unidade Acad\^emica de Matem\'atica\\
   58109-970 Campina Grande, PB - Brazil}\\

{ \scriptsize email: coalves@yahoo.com.br}

\bigskip

{\scriptsize  \bf Jose  V. A. Goncalves }\\
\medskip
 
\scriptsize{ Universidade Federal de Goi\'as\\
   Instituto de Matem\'atica e Estat\'istica\\
   74001-970 Goi\^ania, GO - Brazil}\\

\scriptsize {email: goncalves.jva@gmail.com}
\medskip

\scriptsize{\bf Kaye O. Silva}\\
\scriptsize{ Instituto Federal Goiano\\
N\'ucleo de Estat\'istica, Matem\'atica e Matem\'atica Aplicada \\
74085-010  Uruta\'i, GO - Brazil}
\smallskip

{\scriptsize email:  kayeoliveira@hotmail.com}
\smallskip

\end{flushright}


\begin{thebibliography}{99}
\bibitem{gs-1}  J. V. Goncalves \& C. A.  Santos,  {\it Positive solutions of some quasilinear singular second order equations.}
 Journal of the Australian Mathematical Society  76 (2001) 125-140.

\bibitem{gs-2}  J. V. Goncalves  \&  C. A. Santos,  {\it Classical solutions of singular Monge-Amp\'ere equations in a ball},  Journal of Mathematical Analysis and Applications,   305  (2005)  240-252.

\bibitem{mg}  J. V. Goncalves \& A. L. Melo,   {\it  Multiple sign changing solutions in  a class of quasilinear equations}, Differential Integral Equations,  15   (2002)  147-165.

\bibitem{fn-1} N.  Fukagai  \& K. Narukawa,  {\it  On the existence of multiple positive solutions of quasilinear elliptic eigenvalue problems},   Annali di Matematica    186   (2007) 539-564.

\bibitem{Cheng} Y. Cheng,  {\it On the existence of radial solutions of a nonlinear elliptic equation on
the unit ball}, Nonlinear Anal. 24 (1995) 287 - 307.

\bibitem{JIa} J. Iaia, {\it Radial solutions to a $p$-Laplacian Dirichlet problem}, Applicable Anal. 58 (1995) 335 - 350.


\bibitem{CKr} A. Castro \& A. Kurepa, {\it Infinitely many radially symmetric solutions to a superlinear Dirichlet problem in a ball}, Proc. Amer. Math. Soc. 101 (1987) 57 -64.

\bibitem{CFM} Ph. Cl\'ement, D.G. de Figueiredo \& E. Mitidieri, {\it Quasilinear elliptic equations with critical exponents}, Top. Meth. Nonl. Anal. 7 (1996), 133-170.


\bibitem{CCN} A. Castro, J. Cossio \& J. M. Neuberger, {\it A sign-changing solution for a superlinear Dirichlet problem}, Rocky Mountain J. Math. 27 (4) (1997), 1041-1053.


\bibitem{MR1}M. Mihailescu \& V. Radulescu; {\it Nonhomogeneous Neumann problems in Orlicz-Sobolev spaces,}
C .R. Acad. Sci. Paris, Ser. I 346 (2008), 401-406.

\bibitem{MR2}M. Mihailescu \& V. Radulescu;{\it Existence and multiplicity of solutions for a quasilinear non-
homogeneous problems: An Orlicz-Sobolev space setting,} J. Math. Anal. Appl. 330 (2007),
416-432.

\bibitem{NiSer} W. M. Ni \& J. Serrin, {\it Nonexistence theorems for singular solutions of
quasilinear partial differential equations}, Comm. Pure Appl. Math. 39 (1986) 379 - 399.

\bibitem{Sim} J. Simon, {\it Regularit\'{e} de la solution d'une equation non lineaire  dans ${\r^N}$}, Springer Lecture Notes in Math. \# 665 (Ph. Benilan Editor) (1978), 203 - 227.

\bibitem{St} W. Strauss, {\it Existence of solitary waves in higher dimensions}, Comm. Math. Phys. 55 (1977) 149--162.

\bibitem{Tolk} P. Tolksdorff, {\it On quasilinear boundary value problems
in domains with corners}, Nonlinear Anal. 5 (1981) 721-735.

\bibitem{SWei} R. Saxton \& D. Wei, {\it Radial solutions to a nonlinear $p$-harmonic Dirichlet problem}, Applicable Anal. 51 (1993) 59 - 80.









\end{thebibliography}
\end{document}